\documentclass[11pt]{amsart}
\usepackage{graphicx,color}
\usepackage{natbib}
\newtheorem{thm}{Theorem}[section]
\newtheorem{cor}[thm]{Corollary}
\newtheorem{lem}[thm]{Lemma}

\theoremstyle{definition}

\theoremstyle{remark}
\newtheorem{rem}[thm]{Remark}
\newtheorem{ex}[thm]{Example}
\numberwithin{equation}{section}

\newcommand{\mgi}{M/G/$\infty$ }

\newcommand{\mmone}{M/M/1 }
\newcommand{\mmc}{M/M/C }
\newcommand{\mmi}{M/M/$\infty$ }

\newcommand{\matP}{{\bf P}}
\newcommand{\matQ}{{\bf Q}}
\newcommand{\matR}{{\bf R}}
\newcommand{\matV}{{\bf V}}
\newcommand{\matU}{{\bf U}}
\newcommand{\matD}{{\bf D}}
\newcommand{\matT}{{\bf T}}
\newcommand{\matM}{{\bf M}}
\newcommand{\matW}{{\bf W}}
\newcommand{\matB}{{\bf B}}
\newcommand{\matE}{{\bf E}}
\newcommand{\matI}{{\bf I}}
\newcommand{\matG}{{\bf G}}
\newcommand{\matH}{{\bf H}}
\newcommand{\matPi}{{\bf \Pi}}
\newcommand{\matLambda}{{\bf \Lambda}}

\newcommand{\calN}{{\mathcal N_\Gamma}}
\newcommand{\Po}{\text{Po}}
\newcommand{\dimE}{{K}}
\def\symbtr{\dagger}
\def\transp#1{#1^\symbtr}
\def\ffact#1{^{\underline #1}}
\def\rfact#1{^{\overline  #1}}

\newcommand{\diag}{\text{\bf diag}}

\def\dvec#1(#2){\vbox{\halign{##\cr$\scriptscriptstyle#2\rightarrow$\cr\noalign{\kern1pt\nointerlineskip}$\hfil\displaystyle{#1}\hfil$\cr}}}

\newcommand{\defeq}{:=}
\newcommand{\disteq}{\stackrel{\rm d}{=}}

\newcommand{\R}{\mathbb{R}}
\newcommand{\Rp}{\mathbb{R}^+}
\newcommand{\Rn}{\mathbb{R}^-}
\newcommand{\Rq}{\mathbb{R}^2}
\newcommand{\Z}{\mathbb{Z}}
\newcommand{\E}{\mathbb{E}}

\newcommand{\T}{\mathcal{T}}

\newcommand{\m}[1]{{|#1|}}

\newcommand{\ImF}{n} 
\newcommand{\ImS}{i} 
\newcommand{\ImT}{l} 
\newcommand{\IsF}{k} 
\newcommand{\IsS}{j} 
\newcommand{\IcF}{h} 

\newcommand{\ItlF}{\ImT}              
\newcommand{\ItlS}{\ImF}              
\newcommand{\ItlT}{\ImS}              

\begin{document}

\title[\mmi queues in semi-Markovian RE]{\mmi queues in semi-Markovian random environment}%
\author{B. D'Auria}%
\address{Universidad Carlos III de Madrid, Dpto. de Estadistica, Avda. de la Universidad
30, 28911  Leganés  (Madrid), Spain}%
\email{bernardo.dauria@uc3m.es}%

\thanks{Part of this research took place while the author was still post-doc at EURANDOM, Eindhoven, The Netherlands.}%

\date{}%
\begin{abstract}
In this paper we investigate an \mmi queue whose parameters depend on an external random environment that we assume to be a semi-Markovian process with finite state space. For this model we show a recursive formula that allows to compute all the factorial moments for the number of customers in the system in steady state. The used technique is based on the calculation of the raw moments of the measure of a bidimensional random set. Finally the case when the random environment has only two states is deeper analyzed. We obtain an explicit formula to compute the above mentioned factorial moments when at least one of the two states has sojourn time exponentially distributed.
\end{abstract}
\maketitle
\section{Introduction}
The \mmi queue is one of the simplest model in queueing theory. This is due to the joint situation to have a memory-less arrival process and an infinite set of servers that allows customers to behave independently from each other. This suddenly stops to be true after introducing some correlation between customers. In this paper we achieve that by introducing an independent random environment that modulates the system parameters, i.e. the arrival rate and the server speeds.
Queues with variable service and arrival speeds arise naturally in practice and therefore many classical works can be found. Most of the results deal with the single server queue, see for example \cite{takine:2005}, \cite{ozawa:2004}, \cite{sengupta:1990} and references therein. \cite{neuts:1981} analyzed the  \mmone queue as well as the \mmc queue in random environment by using the matrix-geometric approach while \cite{takine:sengupta:1997} looked at the infinite server queue when only the arrival process was subject to a Markovian modulation.
The infinite server queue in random environment has then been studied by \cite{keilson:servi:1993}, \cite{baykal-gursoy:xiao:2004} and \cite{dauria:2007} in the special case when the random environment is Markovian and has only two states.

In \cite{o'cinneide:purdue:1986} the authors looked at the case when the environment is given by a finite state Markov process and for this case they showed how to compute the factorial moments for the number of customers in the system in steady state. Here we extend their analysis to the case of a semi-Markovian random environment.

This extension is interesting as it makes the model more attractive for application purposes. Indeed, despite its simplicity, the \mmi system is often used to model pure delay systems, such as highways, satellite links or long communication cables, or to approximate the behavior of multi server systems. When these kinds of systems are subject to external influences, such as day time changing rates, it is then helpful to look at extended models, such as the one proposed in this work, in order to analyze or predict their behaviors.

The methodology we use follows the technique developed in \cite{dauria:2007}. It consists on representing the stationary and isolate \mgi system as a Poisson process on $\R^2$ and by computing the number of customers in the system by measuring a deterministic set according to the point process measure (see also \citet[\S 3.3]{resnick:1987} and \cite{dauria:resnick:2006}). In this context the random environment can be expressed as a random modulation of the set and in the special case of exponential distributed service times its measure can be derived by solving a system of stochastic equations, see relation (\ref{eq:A0x.rec}) below. We use this relation to compute the factorial moments for the number of customers in the system at steady state.

\section{Model description}
To start, we define the random environment $\{\Gamma(u),u\in\R\}$ as a semi-Markov chain with values in the finite state space $E=\{1,\ldots,\dimE\}$. We assume that the sojourn time in the state $\IsF\in E$, denoted as $T_\IsF$, is an independent positive random variable whose distribution function has Laplace transform denoted by $\tau_\IsF(s)\defeq\E[e^{-s T_\IsF}]$. In the following we show that the Laplace transform is the only information we need to compute the moments. When the sojourn time in state $\IsF\in E$ expires, the environment jumps to state $\IsS\in E$ with probability $p_{\IsF\IsS}$. Denoting by $\matP\defeq\{p_{\IsF\IsS}\}_{\IsF,\IsS\in E}$ the routing matrix that we assume irreducible and with no loss of generality with $p_{\IsF\IsF}=0$, we can define the reverse routing matrix \begin{equation}\label{eq:def.matQ}
\matQ\defeq \matPi^{-1} \, \transp\matP \, \matPi,
\end{equation}
where $\transp\matP$ denotes the transpose of the matrix $\matP$, $\matPi\defeq\diag(\vec{\pi})$ and $\vec{\pi}$ is the stationary distribution of the Markov chain generated by $\matP$ \citep[see][\S 6.1]{bremaud:1999}.

We assume that when the environment is in state $\IsF\in E$ customers arrive according to a Poisson rate $\lambda_\IsF\ge0$.  Each of them brings an independent request of service, $\sigma$, that is exponentially distributed with rate $\mu$. All servers work at constant speed $\beta_\IsF=\mu_\IsF/\mu\leq1$.
To avoid trivial cases we assume that $\mu,\beta,\lambda>0$ where $\lambda\defeq\max_{\IsF\in E}\lambda_\IsF$ and $\beta\defeq\max_{\IsF\in E}\beta_\IsF$.

By the results in \cite{dauria:2007} the system is stable and we are allowed to study its stationary regime.

We then look at the system at time $0$ and we count the number of customers still in the system. We order them according to their arrival times $\{u_\IcF\}_{\IcF\in\Z}$ with $u_\IcF<u_{\IcF+1}$ and $u_{-1}<0\leq u_0$, and we denote by $G(\sigma) \defeq 1-e^{-\mu\sigma}$, $\sigma>0$, the common exponential distribution function of the $\{\sigma_\IcF\}_{\IcF\in\Z}$.

The $\IcF$-th customer, $\IcF<0$, will be in the system at time $0$ iff its service time, $\sigma_\IcF$, is bigger than the work done by the server it has occupied during the time interval $[u_\IcF,0)$. We denote this quantity by $F_\Gamma(u_\IcF)$ and, as the subscript shows, it is a random quantity that depends on the random environment $\Gamma$. Its value can be computed in the following way,
\begin{equation}\label{eq:def.F_Gamma}
  F_\Gamma(u)\defeq\int_u^0 \beta_{\Gamma(t)} dt, \quad u\leq0.
\end{equation}

Denoting by $N$ the number of customers in stationary regime we have that it is given by
\begin{equation}\label{eq:def.N}
  N=\sum_{\IcF<0} 1\{\sigma_\IcF>F_\Gamma(u_\IcF)\},
\end{equation}
where $1\{\cdot\}$ is the indicator function of the set $\{\cdot\}$.
It is helpful to rewrite the numerable collection of indicator functions appearing in expression (\ref{eq:def.N}) in the following equivalent way
\begin{equation*}
  1\{\sigma_\IcF>F_\Gamma(u_\IcF)\}=\delta_{(u_\IcF,\sigma_\IcF)}(A_\Gamma)
\end{equation*}
where $\delta_{(u,\sigma)}$ denotes a Dirac delta measure with center $(u,\sigma)\in\R\times\Rp$ and the set $A_\Gamma\subset\Rn\times\Rp$ is given by
\begin{equation}\label{eq:def.A_Gamma}
A_\Gamma\defeq\{(u,\sigma)\in\Rn\times\Rp:\sigma>F_\Gamma(u)\}.
\end{equation}
This alternative formulation allows the decoupling of the sequence $\{(u_\IcF,\sigma_\IcF)\}_{\IcF\in\Z}$ and the function $F_\Gamma(u)$ both depending on the realization of the environment $\Gamma$ in the computation of the quantity $N$. Indeed we can express the stationary number of customers in the system in the following way
\begin{equation}\label{eq:def.N.by.calN}
  N=\sum_{\IcF<0} \delta_{(u_\IcF,\sigma_\IcF)}(A_\Gamma)=\calN(A_\Gamma),
\end{equation}
where
\begin{equation}\label{eq:def.calN}
  \calN\defeq\sum_{\IcF\in\Z} \delta_{(u_\IcF,\sigma_\IcF)}
\end{equation}
is a point process  which locates one Dirac delta measure at each arrival point $\{(u_\IcF,\sigma_\IcF)\}_{\IcF\in\Z}$. For the theoretical background and definition of point processes see \citet[][\S 7.1]{daley:vere-jones:1988} or \citet[][\S 3.1]{resnick:1987}. The subscript $\Gamma$ stays to denote that $\calN$ depends on the random environment by the sequence of arrival times $\{u_\IcF\}_{\IcF\in\Z}$.
Indeed given a realization $\gamma$ of the process $\Gamma$, the sequence $(\{u_\IcF\}_{\IcF\in\Z}|\Gamma=\gamma)$ belongs to an inhomogeneous Poisson process with intensity rate $\lambda_{\gamma(u)}$, $u\in\R$. By Proposition 3.8 in \cite{resnick:1987} it follows that $\calN|\Gamma=\gamma$ is still a Poisson process, now on $\R\times\Rp$, with intensity measure
\begin{equation*}
\lambda_\gamma(A)\defeq E[\calN(A)|\Gamma=\gamma]=\int_A \lambda_{\gamma(u)} du \, G(d\sigma), \quad A\subset\R\times\Rp.
\end{equation*}
Finally $\calN$ is a doubly stochastic Poisson process or, more briefly, a Cox process \citep[see][\S8.5]{daley:vere-jones:1988}, i.e. a Poisson process whose intensity measure is itself random and given by
\begin{equation}\label{eq:def.lambda_Gamma}
\lambda_\Gamma(A)\defeq E[\calN(A)|\Gamma]=\int_A \lambda_{\Gamma(u)} du \, G(d\sigma), \quad A\subset\R\times\Rp.
\end{equation}
It is well know that the fidi distributions of a Cox process are of mixed Poisson type \citep[see][Corollary 8.5.II]{daley:vere-jones:1988}, or equivalently that for any set $A\subset\R\times\Rp$
\begin{equation}\label{eq:Cox.process.fidi}
\calN(A)\sim\text{Po}(|A|_\Gamma)
\end{equation}
is a Poisson random variable whose parameter is itself random with value  $|A|_\Gamma=\E[\calN(A)|\Gamma]$. $|A|_\Gamma$ can be geometrically interpreted as the measure of the set $A$ according to the measure $\lambda_\Gamma(\cdot)$, i.e. $|A|_\Gamma=\lambda_\Gamma(A)$.

From relations (\ref{eq:def.N.by.calN}) and (\ref{eq:Cox.process.fidi}) we finally get that
\begin{equation}\label{eq:N}
  N\sim\text{Po}(|A_\Gamma|_\Gamma),
\end{equation}
a mixed Poisson random variable with random parameter $|A_\Gamma|_\Gamma=\lambda_\Gamma(A_\Gamma)$.
\begin{figure}[h]
\begin{center}
    \includegraphics[height=6cm]{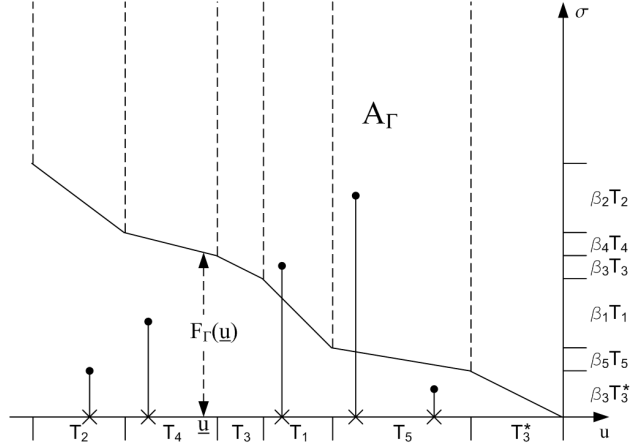}
\caption{Example of realization.}
\label{fg:random.set}
\end{center}
\end{figure}
Figure \ref{fg:random.set} shows an example of realization where the random environment has $K=5$ states: the dots are the centers of the Dirac deltas of the point process $\calN$, while the piecewise linear function $F_\Gamma(u)$ denotes the lower bound of the set of integration $A_\Gamma$. The customers present in the system at time $0$ are then the ones whose dots fall in the set $A_\Gamma$; in the shown example $N=2$.

\begin{ex}
  The easiest case is when the environment process is constant, $K=1$, so that the system reduces to a classical \mmi queue. In this case the set $A_\Gamma$ is deterministic, given by $\{(u,\sigma): u<0,\, \sigma> \beta |u|)\}$. From equation (\ref{eq:def.lambda_Gamma}) we get
  \begin{equation*}
  |A_\Gamma|_\Gamma = \int_{-\infty}^0 \int_{\sigma>\beta |u|} \lambda G(d\sigma) \, du
                    = \frac{\lambda}{\beta} \int_0^\infty 1-G(u) \, du
                    = \frac{\lambda}{\beta} \E[\sigma] = \frac{\lambda}{\beta \mu}
  \end{equation*}
  and we obtain the known information that
  $ N\sim\text{Po}(\frac{\lambda}{\beta \mu})$, i.e. the stationary number of customer in the system is Poisson distributed.
\end{ex}

\section{Computing the factorial moments}
Before beginning to compute the factorial moments of the random variable $N$, it is worthwhile to review some basic results about the different kinds of moments and their relations with the various generating functions. A good reference about the following relations especially in connection with point processes is \cite{daley:vere-jones:1988}, Chapter 5.

Given a random variable $X$, we denote by $\psi_X(s)\defeq\E[e^{s X}]$ its moment generating function and by $\phi_X(z)\defeq\E[z^X]$ its probability generating function.

The factorial moment of order $\ImS$ of $X$, $f_X^{(\ImS)}$ is defined as
$$ f_X^{(\ImS)}\defeq\E\left[X\ffact{\ImS}\right]=\sum_{\ImF=0}^\infty \ImF\ffact{\ImS} \, p_\ImF,$$
where $p_\ImF=\Pr\{X=\ImF\}$ and $\ImF\ffact{\ImS}\defeq \ImF (\ImF-1) \cdots (\ImF-\ImS+1)$ is the falling factorial. It can be directly computed by the $\ImS$-th derivative of the probability generating function, i.e.  $f_X^{(\ImS)}=\lim_{z\to1} \phi_X^{(\ImS)}(z)$.
Knowing the factorial moments of $X$ it is then easy to compute its moments, in the sequel called {\it raw moments} to distinguish them from the factorial ones. Indeed, by taking the expectations on both sides of the following known equivalence \citep{abramowitz:stegun:1964}
$$X^\ImF=\sum_{\ImS=0}^\ImF \mathfrak{S}_\ImF^{(\ImS)} X\ffact{\ImS},$$
where $\mathfrak{S}_\ImF^{(\ImS)}$ is a Stirling Number of the Second Kind, we obtain the following relation between the $\ImF$-th moment of $X$, \hbox{$m_X^{(\ImF)}\defeq\E\left[X^\ImF\right]$} with $m_X^{(0)}\defeq1$, and the factorial moments of order $\ImS\leq \ImF$,
\begin{equation}\label{eq:from.factorial.to.raw.moments}
  m_X^{(\ImF)}=\sum_{\ImS=0}^\ImF \mathfrak{S}_\ImF^{(\ImS)} \, f_X^{(\ImS)}.
\end{equation}

The reverse relation is obtained by using the Stirling Numbers of the First Kind, $\mathfrak{s}_\ImS^{(\ImF)}$
 \citep[see][]{abramowitz:stegun:1964}, that satisfy the following known relation
$$ X\ffact{\ImS}=\sum_{\ImF=0}^\ImS \mathfrak{s}_\ImS^{(\ImF)} X^\ImF,$$
so that, taking the expectations of both sides, finally we get
\begin{equation}\label{eq:from.raw.to.factorial.moments}
  f_X^{(\ImS)}=\sum_{\ImF=0}^\ImS \mathfrak{s}_\ImS^{(\ImF)} \, m_X^{(\ImF)}.
\end{equation}

It is interesting to notice that relation (\ref{eq:from.factorial.to.raw.moments}) comes directly from using the fact that  $\psi_X(s)=\phi_X (e^s)$ and that $m_X^{(\ImF)}=\lim_{s\to0} \psi_X^{(\ImF)}(s)$.
Indeed,
$$\lim_{s\to0} \psi_X^{(\ImF)}(s)=\lim_{s\to0} \frac {d^\ImF} {d s^\ImF} \phi_X (e^s)=\sum_{\ImS=0}^\ImF  \mathfrak{S}_\ImF^{(\ImS)} \, \phi_X^{(\ImS)}(1),$$
where in the last equation we used {Fa\'a di Bruno's formula} for the expansion of derivatives of order $\ImF$ for composition of functions \citep[see][]{abramowitz:stegun:1964} and the fact that $\lim_{s\to0} \frac {d^\ImF} {d s^\ImF} e^{s} = 1$.

A random variable $X$ is called mixed Poisson when there exists a\break non-negative random variable $Y$ such that $X\disteq\Po(Y)$, or equivalently $\left(X|Y=y\right)\disteq\Po(y)$, where the operator $\disteq$ denotes equality in distribution. In the case $X$ were a mixed Poisson random variable we would have that
$$\phi_X(z)=\psi_Y(z-1),$$
so that taking the derivatives of order $\ImF$, we get
$$\lim_{z\to1} \phi_X^{(\ImF)}(z)=\lim_{z\to1} \psi_Y^{(\ImF)}(z-1)=\lim_{s\to0} \psi_Y^{(\ImF)}(s),$$
or, in other words, that the factorial moments of $X$ are directly the raw moments of $Y$,
$$f_X^{(\ImF)}= m_Y^{(\ImF)},$$
and the latter often are easier to compute.

This is exactly what happens in our case where, as shown by relation (\ref{eq:N}), $N$ is a mixed Poisson and that is why we are interested into its factorial moments rather then directly its raw moments. Indeed we have that the following relation holds
\begin{equation}\label{eq:fm.N.rm.A}
  f_N^{(\ImF)}= m_{|A_\Gamma|_\Gamma}^{(\ImF)},
\end{equation}
and our task reduces to the computation of the raw moments of the measure of the random set $A_\Gamma$.

\section{Computing the raw moments of $|A_\Gamma|_\Gamma$}

In this section we compute the raw moments of the measure of the set $A_\Gamma$, defined in (\ref{eq:def.A_Gamma}), when measured by the random intensity measure $\lambda_\Gamma$, defined in (\ref{eq:def.lambda_Gamma}).
We use a fixed point technique and to this aim  we look at a modified environment process, $\Gamma_0$, that is the Palm version of the process $\Gamma$, i.e. we assume that at time $0$ it has a transition. We denote by $\IsF\in E$ the last state it has assumed before $0$, i.e. $\IsF\defeq\Gamma_0(0^-)$,  and by $T_\IsF$ its corresponding sojourn time.
While, as depicted in Figure \ref{fg:random.set}, for the process $\Gamma$ the sojourn time in the last state before $0$ would be given by a residual sojourn time, for the process $\Gamma_0$ it is distributed as any other sojourn time corresponding to the same state.
We define by $A_{0\IsF}\defeq(A_{\Gamma_0}|\Gamma_0(0^-)=\IsF)$, $\IsF\in E$, the set $A_{\Gamma_0}$ conditioned to the event that the last state occupied by the environment before $0$ is the state $\IsF$, and we call $\m{A_{0\IsF}}$ its measure, i.e. $\m{A_{0\IsF}}\defeq(\lambda_{\Gamma_0}(A_{\Gamma_0})|\Gamma_0(0^-)=\IsF)$.

\begin{figure}[h]
\begin{center}
    \includegraphics[height=6cm]{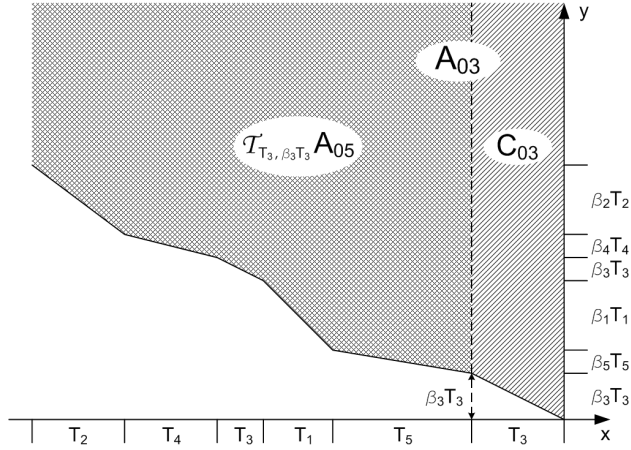}
\caption{Decomposition of $A_{03}$ as $C_{03} \cup \T_{T_3,\beta_3 T_3} A_{05}$.}
\label{fg:decomposition}
\end{center}
\end{figure}

Figure \ref{fg:decomposition} shows an example of the set $A_{0\IsF}$ when $\IsF=3$, together with its decomposition in the set $C_{0\IsF}$ and the set $\T_{T_3,\beta_3 T_3} A_{05}$. To this we have defined by $C_{0\IsF}$ the restriction of the set $A_{0\IsF}$ up the last transition of the process $\Gamma_0$ before time $0$, i.e.
\begin{align}\label{eq:def.C0x}
\phantom{|}C_{0\IsF}\phantom{|}&\defeq A_{0\IsF} \cap \{(x,y)\in\Rn\times\Rp| |x|<T_\IsF\},\\ \m{C_{0\IsF}}&\defeq(\lambda_{\Gamma_0}(C_{0\IsF})|\Gamma_0(0^-)=\IsF)\nonumber
\end{align}
and by $\T_{s,t} A$ the $(-s, t)$-translated version of the set $A$, i.e.
\begin{equation}\label{eq:def.Tst}
\T_{s,t} A \defeq \{(x,y)\in\Rq| (x+s,y-t)\in A\}.
\end{equation}
We denote by $\IsS\in E$ the state of the environment before the last transition before time $0$, i.e. $\IsS\defeq\Gamma_0(-T_\IsF^-)$, so that, $-T_\IsF$ being a regeneration point for the process $\Gamma_0$, we have the independence of the sets $C_{0\IsF}$ and $\T_{T_\IsF, \beta_\IsF T_\IsF} A_{0\IsS}$ conditioned to the values of the states before and after the transitions, i.e. $\IsS$ and $\IsF$. $\beta_\IsF T_\IsF=F_{\Gamma_0}(-T_\IsF)$ is the exact amount of work the non-empty servers have done during the time interval $[-T_\IsF,0)$ being in state $\IsF$.

By noticing that the set $\T_{T_\IsF, \beta_\IsF T_\IsF} A_{0\IsS}$ has measure equal in distribution to $\m{\T_{0, \beta_\IsF T_\IsF} A_{0\IsS}}\defeq(\lambda_{\Gamma_0}(\T_{0, \beta_\IsF T_\IsF} A_{0\IsS})|\Gamma_0(0^-)=\IsS)$, we can write down the following set of stochastic equations
\begin{equation}\label{eq:A0x}
 \m{A_{0\IsF}} \disteq \m{C_{0\IsF}} + \sum_{\IsS=1}^\dimE 1\{\IsF\leftarrow \IsS\} \m{\T_{0,\beta_\IsF T_\IsF} A_{0\IsS}},
\end{equation}
where the indicator function $1\{\IsF\leftarrow \IsS\}$ selects the backward state transition of the environment from the state $\IsF$ to the state $\IsS$; this would happen, according to definition (\ref{eq:def.matQ}), with probability $q_{\IsF\IsS}$.

Thanks to the fact that along the vertical axis the measure function is given by $G$ that is exponential we have that the following result holds:
\begin{lem}\label{lm:Tst}
  Given the transformation $\T_{s,t}$, defined in (\ref{eq:def.Tst}), we have that
  \begin{equation}
    |\T_{0,t} A_\Gamma|_\Gamma= e^{-\mu t} |A_\Gamma|_\Gamma,
  \end{equation}
  for any random set $A_\Gamma:\Gamma \to \mathcal{B}(\Rn\times\Rp)$.
\end{lem}
\begin{proof}
By using definition (\ref{eq:def.lambda_Gamma}) we have
\begin{align*}
  \lambda_{\Gamma}(\T_{0,t} A_\Gamma)&=\int_{\T_{0,t} A_\Gamma} \lambda_{\Gamma(u)} du \, e^{-\mu\sigma}d\sigma\\
  &=e^{-\mu t} \int_{\T_{0,t} A_\Gamma} \lambda_{\Gamma(u)} du \,   e^{-\mu(\sigma-t)}d\sigma=e^{-\mu t} \lambda_{\Gamma}(A_\Gamma).
\end{align*}
\end{proof}

By using Lemma \ref{lm:Tst}, equation (\ref{eq:A0x}) simplifies in the following
\begin{equation}\label{eq:A0x.rec}
 \m{A_{0\IsF}} \disteq \m{C_{0\IsF}} + e^{-\mu \beta_\IsF T_\IsF} \sum_{\IsS=1}^\dimE 1\{\IsF\leftarrow \IsS\} \m{A_{0\IsS}},
\end{equation}
that is the starting point to prove the following main result:
\begin{thm}
  Let us define $\vec m_0^{(\ImF)}\in\R^\dimE$ as the column vector whose $\IsF$-th coordinate is the $\ImF$-th moment of the random variable $\m{A_{0\IsF}}$, i.e. \hbox{$m_{0\IsF}^{(\ImF)}\defeq m^{(\ImF)}_{\m{A_{0\IsF}}}$} then the following relation holds
  \begin{equation}\label{eq:main.result.imp}
    \sum_{\ImS=0}^\ImF (-1)^\ImS \binom{\ImF}{\ImS} \matR^{\ImF-\ImS} \matB_\ImF \, \vec m_0^{(\ImS)}=0,
  \end{equation}
  where $\matR\defeq\diag\left(\rho_\IsF\right)$, $\rho_\IsF\defeq\lambda_\IsF/\mu_\IsF$ and the matrix \hbox{$\matB_\ImF\defeq\diag\left(\tau_\IsF^{-1}(\ImF\mu_\IsF)\right) - \matQ$}. The matrix $\matB_\ImF$, $\ImF>0$, is invertible and therefore it is possible to express the $\ImF$-th moment vector $\vec m_0^{(\ImF)}$ in terms of the previous ones, $\vec m_0^{(\ImS)}$, $\ImS=0,\ldots,\ImF-1$, in the following way
  \begin{equation}\label{eq:main.result.exp}
    \vec m_0^{(\ImF)}= \sum_{\ImS=0}^{\ImF-1} (-1)^{\ImF-1-\ImS} \binom{\ImF}{\ImS} \matB_\ImF^{-1} \matR^{\ImF-\ImS} \matB_\ImF \, \vec m_0^{(\ImS)}.
  \end{equation}
\end{thm}
\begin{proof}
  We first compute the values of the variable $\m{C_{0\IsF}}$ in the following way
  $$\m{C_{0\IsF}}= \lambda_\IsF \int_0^{T_\IsF} e^{-\mu_\IsF x} dx = \rho_\IsF (1- e^{-\mu_\IsF T_\IsF}).$$
  Then substituting its value in equation (\ref{eq:A0x.rec}), it gives
\begin{equation}
  \m{A_{0\IsF}} \disteq \rho_\IsF (1- e^{-\mu_\IsF T_\IsF})+
e^{-\mu_\IsF T_\IsF} \sum_{\IsS=1}^\dimE 1\{\IsF \leftarrow \IsS\} \m{A_{0\IsS}},
\end{equation}
that can be rewritten as
\begin{equation}\label{eq:A.0x}
 \m{A_{0\IsF}}-\rho_\IsF \disteq  \sum_{\IsS=1}^\dimE 1\{\IsF \leftarrow \IsS\} (\m{A_{0\IsS}} -\rho_\IsF) e^{-\mu_\IsF T_\IsF}.
\end{equation}
We denote by $\psi_{0\IsF}(s)\defeq\E\left[e^{s \m{A_{0\IsF}}}\right]$ the moment generating function of $\m{A_{0\IsF}}$ so that applying the exponential function to both members of equation (\ref{eq:A.0x}) previously multiplied by $s$ and then taking the expectation, we obtain
\begin{align*}
  \psi_{0\IsF}(s) e^{-s \rho_\IsF}
  &= \E\left[\sum_{\IsS=1}^\dimE q_{\IsF\IsS} e^{s(\m{A_{0\IsS}} -\rho_\IsF) e^{-\mu_\IsF T_\IsF}}\right]\\
  &= \E\left[\sum_{\IsS=1}^\dimE q_{\IsF\IsS} \psi_{0\IsS}(s e^{-\mu_\IsF T_\IsF}) e^{-s \rho_\IsF e^{-\mu_\IsF T_\IsF}}\right].
\end{align*}
Last expression can be written in matrix form in the following way
\begin{equation}\label{eq:A_0x}
 e^{-s \matR} \vec\psi_0(s) =\E\left[e^{-s \matR \matT} (\matQ \vec\psi_0)(s\matT)\right],
\end{equation}
where $\matT\defeq\diag(e^{-\mu_\IsF T_\IsF})$ and where with notation $\vec v(\matW)$, with $\matW$ a diagonal matrix, we denote a vector whose $\IsF$-th component is $v_\IsF(w_{\IsF\IsF})$.
We use then the following matrix formulas for derivatives
\begin{equation}
  D^{(\ImF)}[e^{-s \matW}\vec v(s)]=\sum_{\ImS=0}^\ImF (-1)^{\ImF-\ImS} \binom{\ImF}{\ImS} e^{-s\matW} \matW^{\ImF-\ImS} D^{(\ImS)}[{\vec v}(s)],
\end{equation}
and
\begin{equation}
  D^{(\ImF)}[\vec v(s \matW)]=\matW^\ImF {\vec v}^{\,(\ImF)}(s \matW),
\end{equation}
to compute the $\ImF$-th derivative of both sides of equation (\ref{eq:A_0x}) so that
\begin{align*}
\sum_{\ImS=0}^\ImF (-1)^{\ImF-\ImS} \binom{\ImF}{\ImS} &e^{-s\matR} \matR^{\ImF-\ImS} \vec \psi_0^{\,(\ImS)}(s)=\\
&=\E\left[\sum_{\ImS=0}^\ImF (-1)^{\ImF-\ImS} \binom{\ImF}{\ImS} e^{-s\matR\matT} \matR^{\ImF-\ImS} \matT^{\ImF-\ImS} D^{(\ImS)}[{\matQ \vec\psi_0}(s \matT)]\right]\\
&=\E\left[\sum_{\ImS=0}^\ImF (-1)^{\ImF-\ImS} \binom{\ImF}{\ImS} e^{-s\matR\matT} \matR^{\ImF-\ImS} \matT^{\ImF} (\matQ \vec\psi_0^{\,(\ImS)})(s \matT)\right].
\end{align*}
Remembering that \hbox{$ \vec m_0^{(\ImF)}=\lim_{s\to0} \vec \psi_0^{\,(\ImF)}(s)$} and taking the limit of last expression as $s\to0$, we get
\begin{equation}
\sum_{\ImS=0}^\ImF (-1)^{\ImF-\ImS} \binom{\ImF}{\ImS} \matR^{\ImF-\ImS} \vec m_0^{(\ImS)}
=\E\left[\sum_{\ImS=0}^\ImF (-1)^{\ImF-\ImS} \binom{\ImF}{\ImS} \matR^{\ImF-\ImS} \matT^\ImF \matQ \, \vec m_0^{(\ImS)}\right].
\end{equation}
Multiplying on the left side by $(-1)^{-\ImF}\E[\matT^{\ImF}]^{-1}$, the last expression can be easily rearranged in
\begin{equation}\label{eq:explicit.Bx}
\sum_{\ImS=0}^\ImF (-1)^{-\ImS} \binom{\ImF}{\ImS} \matR^{\ImF-\ImS} [\E\left[\matT^\ImF\right]^{-1}-\matQ ] \vec m_0^{(\ImS)}=0,
\end{equation}
that gives the result. The invertibility of the matrix $\matB_\ImF$ for $\ImF>0$ comes from Lemma \ref{lm:invertibility}.
\end{proof}
It is remarkable that it is possible to express equation (\ref{eq:main.result.imp}) in terms of the forward transition chain $\matP$. The result is contained in the following corollary whose proof comes from simple matrix computations.
\begin{cor}
  A result similar to equation (\ref{eq:main.result.imp}) is valid for the row vector $ \vec m^{\symbtr(\ImF)}_0 \defeq \transp{\left({\vec m}_0^{(\ImF)}\right)}$, that involves the matrix $\matP$ instead of the matrix $\matQ$, i.e.
  \begin{equation}\label{eq:main.result.imp.P}
    \sum_{\ImS=0}^\ImF (-1)^\ImS \binom{\ImF}{\ImS}  \vec m^{\symbtr(\ImS)}_0 \, \matPi \,
    \matB'_\ImF \, \matR^{\ImF-\ImS} =0,
  \end{equation}
  where $\matB'_\ImF\defeq\diag\left(\tau_\IsF^{-1}(\ImF\mu_\IsF)\right) - \matP$.
  The matrix $\matB'_\ImF$ is non-singular when $\ImF>0$.
  \end{cor}

Given the raw moments of the $\m{A_{\Gamma_0}}$, we can successively compute the moments of the measures of the sets $A_\IsF\defeq (A_\Gamma|\Gamma(0)=\IsF)$, $\IsF\in E$. Following previous definitions we define $m^{(\ImF)}_\IsF\defeq m^{(\ImF)}_{\m{A_\IsF}}$. Similarly to equation (\ref{eq:A0x}) we have the following equation
\begin{equation}\label{eq:Ax}
 \m{A_{\IsF}} \disteq \m{C^*_\IsF} + \sum_{\IsS=1}^\dimE 1\{\IsF\leftarrow \IsS\} \m{\T_{0,\beta_\IsF T^*_\IsF} A_{0\IsS}},
\end{equation}
with $\m{C^*_\IsF}=\rho_\IsF(1-e^{-\mu_\IsF T^*_\IsF})$. $T^*_\IsF$ refers to a residual sojourn time of the environment in state $\IsF\in E$. We define by $\tau^*_\IsF(s)$ the Laplace transform of the distribution function of $T^*_\IsF$ and it is related to the one of $T_\IsF$, $\tau_\IsF(s)$, by the relation $\tau^*_\IsF(s)=\bar \tau_\IsF (1-\tau_\IsF(s))/s$, with $\bar \tau_\IsF\defeq\E[T_\IsF]^{-1}$. For the vector of raw moments $\vec m^{(\ImF)}$ the following theorem holds.
\begin{thm}
  The vector $\vec m^{(\ImF)}$ satisfies the following relation with the vector  $\vec m_0^{(\ImF)}$
\begin{equation}\label{eq:result.imp}
\sum_{\ImS=0}^\ImF (-1)^{\ImS} \binom{\ImF}{\ImS} \matR^{\ImF-\ImS} [\vec m^{(\ImS)}-\matE_\ImF \vec m_0^{(\ImS)}]=0,
\end{equation}
with $\matE_\ImF\defeq\diag(\tau^*_\IsF(\ImF\mu_\IsF)/\tau_\IsF(\ImF\mu_\IsF))$. Therefore the vector $\vec m^{(\ImF)}$ can be computed from the previous moments $\{\vec m^{(\ImS)}\}_{\ImS<\ImF}$ and the corresponding vectors $\{\vec m_0^{(\ImS)}\}_{\ImS\leq \ImF}$ in the following way
\begin{equation}\label{eq:result.exp}
\vec m^{(\ImF)}=\matE_\ImF \vec m_0^{(\ImF)} +\sum_{\ImS=0}^{\ImF-1} (-1)^{\ImF-1-\ImS} \binom{\ImF}{\ImS} \matR^{\ImF-\ImS} [\vec m^{(\ImS)}-\matE_\ImF \vec m_0^{(\ImS)}],
\end{equation}
finally
\begin{equation}\label{eq:mom.A.Gamma}
    m^{(\ImF)}\defeq m^{(\ImF)}_{A_\Gamma}=\vec  m^{(\ImF)} \vec \pi.
\end{equation}
\end{thm}
\begin{proof}
Starting by equation (\ref{eq:Ax}) and following the same calculations that brought us from equation (\ref{eq:A0x}) to equation (\ref{eq:explicit.Bx}), we get
\begin{equation}
\sum_{\ImS=0}^\ImF (-1)^{\ImS} \binom{\ImF}{\ImS} \matR^{\ImF-\ImS} [\E\left[{\matT^*}^\ImF\right]^{-1} \vec m^{(\ImS)}-\matQ  \vec m_0^{(\ImS)}]=0,
\end{equation}
that after subtracting  equation (\ref{eq:explicit.Bx}) gives
\begin{equation}
\sum_{\ImS=0}^\ImF (-1)^{\ImS} \binom{\ImF}{\ImS} \matR^{\ImF-\ImS} [\E\left[{\matT^*}^\ImF\right]^{-1} \vec m^{(\ImS)}-\E\left[{\matT}^\ImF\right]^{-1}  \vec m_0^{(\ImS)}]=0,
\end{equation}
and by multiplying on the left by $\E\left[{\matT^*}^\ImF\right]$ we finally get the result.
\end{proof}

In order to check our results we compare equation (\ref{eq:main.result.imp.P}) for the exponential case with results in \cite{o'cinneide:purdue:1986} here repeated in formula (\ref{eq:o'cinneide:purdue:1986:result}). For this case since $T^*_\IsF \disteq T_\IsF$, we have that $\vec m^{(\ImF)} = \vec m_0^{(\ImF)}$.
\begin{rem}
  It is worth to notice that in \cite{o'cinneide:purdue:1986}, they actually computed the factorial moments of the random row vector \hbox{$(N\, 1\{\Gamma_0=\IsF\})_{\IsF\in E}$} while here we compute the factorial moments of the row vector \hbox{$(N | \Gamma_0=\IsF)_{\IsF\in E}$}. This explains the presence, in formula (\ref{eq:o'cinneide:purdue:1986:result}), of the additional factor given by matrix $\matPi$.
\end{rem}
\begin{cor}
  In case the sojourn times $T_\IsF$ are exponentially distributed with parameters $\bar \tau_\IsF$ we have that with $\ImF>0$
  \begin{equation}\label{eq:o'cinneide:purdue:1986:result}
    (\vec m^{\symbtr(\ImF)} \matPi)\, (\ImF \matM-\matG)= \ImF \, (\vec m^{\symbtr(\ImF-1)} \matPi) \, \matLambda
  \end{equation}
  where $\matM\defeq\diag(\mu_\IsF)$, $\matLambda\defeq\diag(\lambda_\IsF)$ and $\matG\defeq\bar\matT (\matP-\matI)$,  with \hbox{$\bar\matT\defeq\diag(\bar \tau_\IsF)$}, is the generator of the Markovian Environment.
\end{cor}
\begin{proof}
  When the sojourn times  $T_\IsF$ are exponentially distributed we have that
  $$\E[\matT^{\ImF}]^{-1}=\ImF \matM \bar\matT^{-1}+\matI,$$
  that implies
  $$\bar\matT\matB_\ImF=\ImF \matM-\matH,$$
  with $\matH\defeq\bar\matT (\matQ-\matI)$ being the generator of the reverse-time Markov process.
  By multiplying on the left both sides of equation (\ref{eq:main.result.imp}) by $\bar\matT$ and noticing that it can commute with the powers of the matrix $\matR$, we can rewrite (\ref{eq:main.result.imp}) as
  \begin{equation}\label{eq:to.prove}
      \sum_{\ImS=0}^{\ImF} (-1)^{\ImF-\ImS} \binom{\ImF}{\ImS} \matR^{\ImF-\ImS} [\ImF \matM-\matH] \, \vec m^{(\ImS)}=0.
  \end{equation}

    By Lemma \ref{lm:uniq.sol} with $\matU=\matM$, $\matD=\matR$, $\matV_\ImF=\matH$ and $\vec v^{\,(\ImF)}=\vec m^{(\ImF)}$, and imposing $\vec m^{(0)}=\vec1$, we get that the unique solution of (\ref{eq:to.prove}) is given by
   \begin{equation*}
   (\ImF \matM-H) \vec m^{(\ImF)}= \ImF \,\matR \matM \, \vec m^{(\ImF-1)},
   \end{equation*}
  that transposed reduces to
  \begin{equation*}
     \vec m^{\symbtr(\ImF)} (\matPi\,\matPi^{-1})(\ImF \matM-\transp\matH)= \ImF \, \vec m^{\symbtr(\ImF-1)} \matLambda.
  \end{equation*}
  Multiplying on the right side by $\matPi$ and simplifying we get
  \begin{equation*}
     \vec m^{\symbtr(\ImF)}\matPi(\ImF \matM- \matPi^{-1}\transp\matH\matPi)= \ImF \, \vec m^{\symbtr(\ImF-1)} \matPi \matLambda,
  \end{equation*}
  that gives the result after noticing that $\matG=\matPi^{-1}\transp\matH\matPi$.
\end{proof}
\section{Some explicit formulas - Case $\dimE=2$}
Formulas (\ref{eq:main.result.exp}) and (\ref{eq:result.exp}) show that generally to find the $\ImF$-th moment of the random number of users in the system involves in a complex way the knowledge of all previous moments. Reversely the exponential case, that was already solved in \cite{o'cinneide:purdue:1986}, is easier as the $\ImF$-th vector of moments is related only by a factor to the $(\ImF-1)$-th one. That was anyway hidden in a non-trivial way in formula (\ref{eq:main.result.exp}) so that there could be some other special cases where an easier expression holds.

In this section we have a look to the case when the environment has only two stages, i.e. $\dimE=2$.

This is a very special case and when the sojourn times are all assumed exponentially distributed, it is known how to compute the complete distribution of the number of customers in the system at steady state (see \cite{keilson:servi:1993}, \cite{baykal-gursoy:xiao:2004} and \cite{dauria:2007}).

We give for this case an explicit formula to calculate the factorial moments in terms of the Laplace transform of the sojourn time in state $1$, when the sojourn time in state $2$ is exponentially distributed.

By rewriting in more explicit form equation (\ref{eq:A.0x}) for the case $\dimE=2$ we get
\begin{align}
 \m{A_{01}}-\rho_1 &\disteq (\m{A_{02}} -\rho_1) e^{-\mu_1 T_1}\label{eq:A0x.case2.1}\\
 \m{A_{02}}-\rho_2 &\disteq (\m{A_{01}} -\rho_2) e^{-\mu_2 T_2}.\label{eq:A0x.case2.2}
\end{align}
We define $\tilde m_{0\IsF}^{(\ImF)}\defeq\E[(\m{A_{0\IsF}}-\rho_1)^\ImF]$ and take the mean of the $\ImF$-powers of expression (\ref{eq:A0x.case2.1}) so getting
\begin{equation}\label{eq:mom1.case2}
 \tilde m_{01}^{(\ImF)}= \tilde m_{02}^{(\ImF)} \tau_1(\ImF \mu_1).
\end{equation}
By adding and subtracting $\rho_1$ to both sides of equation (\ref{eq:A0x.case2.2}) we get
\begin{equation*}
( \m{A_{02}}-\rho_1)-\rho_*\disteq (\m{A_{01}}-\rho_1)e^{-\mu_2 T_2} -\rho_*  e^{-\mu_2 T_2},
\end{equation*}
with $\rho_*=\rho_2-\rho_1$. Then using equation (\ref{eq:A0x.case2.1}) we obtain a recursive equation involving only $\m{A_{02}}-\rho_1$,
\begin{equation*}
( \m{A_{02}}-\rho_1)-\rho_*\disteq (\m{A_{02}} -\rho_1) e^{-\mu_1 T_1} e^{-\mu_2 T_2} -\rho_*  e^{-\mu_2 T_2}.
\end{equation*}
Taking the $\ImF$-th power and then the expectation of both sides we get
\begin{align*}
 \sum_{\ImS=0}^\ImF (-1)^{\ImF-\ImS} \binom{\ImF}{\ImS} \rho_*^{\ImF-\ImS} \, \tilde m_{02}^{(\ImS)}
 &= \tau_2(\ImF\mu_2) \, \E[\big((\m{A_{02}} -\rho_1) e^{-\mu_1 T_1}-\rho_*\big)^\ImF] \\
 &= \tau_2(\ImF\mu_2) \sum_{\ImS=0}^\ImF (-1)^{\ImF-\ImS} \binom{\ImF}{\ImS} \rho_*^{\ImF-\ImS} \tau_1(\ImS\mu_1) \, \tilde m_{02}^{(\ImS)},
\end{align*}
that, taking into account equation (\ref{eq:mom1.case2}), can be rearranged to get the following
\begin{equation}\label{eq:case2.rec.form}
 \sum_{\ImS=0}^\ImF (-1)^{\ImF-\ImS} \binom{\ImF}{\ImS} \rho_*^{\ImF-\ImS} (\tau_2^{-1}(\ImF\mu_2) \tau_1^{-1}(\ImS\mu_1)-1)\, \tilde m_{01}^{(\ImS)}=0.
\end{equation}
\begin{thm}
  Assuming that the sojourn times of state $2$ are exponentially distributed, i.e. $T_2\sim\text{Exp}(\bar \tau_2)$, the solution of formula (\ref{eq:case2.rec.form}) is given by
  \begin{equation}\label{eq:tilde.m01}
    \tilde m_{01}^{(\ImF)}= \left(\frac{\mu_2 \, \rho_* }{\bar \tau_2}\right)^\ImF  \, \prod_{\ImS=1}^\ImF \frac{\ImS \, \tau_1^{-1}((\ImS-1)\mu_1)}{\tau_1^{-1}(\ImS\mu_1)\tau_2^{-1}(\ImS\mu_2)-1},
  \end{equation}
  and therefore
  \begin{equation}\label{eq:tilde.m02}
    \tilde m_{02}^{(\ImF)}=\left(\frac{\mu_2 \, \rho_* }{\bar \tau_2}\right)^\ImF  \tau_1^{-1}(\ImF \mu_1)\, \prod_{\ImS=1}^\ImF \frac{\ImS \, \tau_1^{-1}((\ImS-1)\mu_1)}{\tau_1^{-1}(\ImS\mu_1)\tau_2^{-1}(\ImS\mu_2)-1}.
  \end{equation}
  Finally
  \begin{equation}\label{eq:m0x}
    m_{0\IsF}^{(\ImF)}=\sum_{\ImS=0}^\ImF \binom{\ImF}{\ImS} \rho_1^{\ImF-\ImS} \tilde m_{0\IsF}^{(\ImS)}, \quad k=1,\,2.
  \end{equation}
\end{thm}
\begin{proof}
  Substituting $\tau_2^{-1}(s)=1+s/{\bar \tau_2}$ in equation (\ref{eq:case2.rec.form}) and rearranging it, we get
  \begin{equation}
 \sum_{\ImS=0}^\ImF (-1)^{\ImF-\ImS} \binom{\ImF}{\ImS} \rho_*^{\ImF-\ImS} \left(\ImF \frac{\mu_2}{\bar \tau_2} -(\tau_1(\ImS\mu_1)-1)\right)\frac{\tilde m_{01}^{(\ImS)}}{\tau_1(\ImS\mu_1)}=0.
\end{equation}
By applying Lemma \ref{lm:uniq.sol} in the scalar case, with $\matU=(\mu_2/{\bar \tau_2})$, $\matD=(\rho_*)$, $\matV_\ImF=(\tau_1(\ImF\mu_1)-1)$ and $\vec v^{\,(\ImF)}=(\tilde m_{01}^{(\ImF)}/\tau_1(\ImF\mu_1))$, we notice that a set of solutions is given by
  \begin{equation}\label{eq:lem.hyp}
     \left(\tau_2^{-1}(\ImF \mu_2)-\tau_1(\ImF\mu_1)\right) \frac{\tilde m_{01}^{(\ImF)}}{\tau_1(\ImF\mu_1)}= \ImF \,\rho_* \frac{\mu_2}{\bar \tau_2} \, \frac{ \tilde m_{01}^{(\ImF-1)}}{\tau_1((\ImF-1)\mu_1)},
  \end{equation}
that is then uniquely defined given that $\tilde m_{01}^{(0)}=1$. Therefore equation (\ref{eq:tilde.m01}) holds. Equation (\ref{eq:tilde.m02}) results by applying (\ref{eq:mom1.case2}) to (\ref{eq:tilde.m01}) and finally equation (\ref{eq:m0x}) comes from the fact that $m_{0\IsF}^{(\ImF)}=\E[((\m{A_{0\IsF}}-\rho_1)+\rho_1)^\ImF]$.
\end{proof}
\begin{ex}\label{ex:known.result}
{\bf  Case $T_1\sim\text{Exp}(\bar\tau_1)$.} In this special case equation (\ref{eq:tilde.m01}) simplifies in
$$\tilde m_{01}^{(\ImF)}= \rho_*^\ImF \frac{(\bar\tau_1/\mu_1)\rfact{\ImF}}{(\bar\tau_1/\mu_1+\bar\tau_2/\mu_2+1)\rfact{\ImF}},$$
with $\ImS\rfact{\ImF}\defeq \ImS (\ImS+1) \cdots (\ImS+\ImF-1)$ being the rising factorial (in \cite{abramowitz:stegun:1964} it is denoted by $(\ImS)_\ImF$). Therefore the moment generating function of $\m{A_{01}}-\rho_1$ is given by the Kummer function $M(\bar\tau_1/\mu_1,\, \bar\tau_1/\mu_1+\bar\tau_2/\mu_2+1,\, \rho_* s)$ (see \cite{abramowitz:stegun:1964}), in accordance to what is shown in \cite{baykal-gursoy:xiao:2004} and in \cite{dauria:2007} [in there, it is denoted by $\phi_{\rm ON}(-s)$].
\end{ex}
The following example is a new result that generalizes the one of Example \ref{ex:known.result}.
\begin{ex}
{\bf  Case $T_1\sim\text{Gamma}(\kappa,\bar\tau_1^{\, -1})$.} For this case we have that $\tau_1^{-1}(k\mu_1)=(1+k \mu_1/{\bar\tau_1})^\kappa$. Therefore equation (\ref{eq:tilde.m01}) simplifies in
$$\tilde m_{01}^{(\ImF)}=  \frac{\ImF! \rho_*^\ImF [(\bar\tau_1/\mu_1)\rfact{\ImF}]^\kappa}{\prod_{\ImS=1}^\ImF \left[(\bar\tau_1/\mu_1+\ImS)^\kappa(\bar\tau_2/\mu_2+\ImS)-(\bar\tau_1/\mu_1)^\kappa(\bar\tau_2/\mu_2)\right]}.$$
\end{ex}
\section{Conclusions}
In this paper we showed that using a matrix-geometric approach it is possible to solve the problem to find the factorial moments of the random number of customers in an \mmi system when its parameters are modulated by a semi-Markovian random environment. We showed that this is possible by looking at this random variable as the random measure of a bidimensional random set by a mixed Poisson process. Finally the case when the environment has only 2 states is more deeply investigated and it is shown that explicit formulas are obtainable given that one state has exponential sojourn times. It is then plausible to believe that for this last case it would be possible to get an explicit expression for the complete characteristic function.
\appendix{}
\section{Technical Lemmas}
\begin{lem}\label{lm:uniq.sol}
   Given the matrices  $\matU,\matD,\{\matV_\ItlS\}_{\ItlS\geq0}$ such that for any \hbox{$0<\ItlS\leq\ItlF$} the matrix $(\ItlS \matU-\matV_\ItlS)$ is invertible than the system of equations
  \begin{equation}\label{eq:lm.result}
  \sum_{\ItlT=0}^{\ItlS}(-1)^{\ItlS-\ItlT} \binom{\ItlS}{\ItlT} \matD^{\ItlS-\ItlT} [\ItlS \matU-\matV_\ItlT] \vec v^{\,(\ItlT)}= 0, \quad 0<\ItlS\leq\ItlF
  \end{equation}
   has a family of solutions $\vec v^{\,(\ItlS)}$, $0<\ItlS\leq\ItlF$, given by
  \begin{equation}\label{eq:lem.hyp}
     (\ItlS \matU-\matV_\ItlS) \vec v^{\,(\ItlS)}= \ItlS \,\matD \matU \, \vec v^{\,(\ItlS-1)},
  \end{equation}
  that reduces to a unique solution once given the vector $\vec v^{\,(0)}$ that has to satisfy the relation $\matV_0 \vec v^{\,(0)}=0$.
\end{lem}
\begin{proof}
  The proof of the lemma is immediate once we prove that a set of vectors obeying to  relations (\ref{eq:lem.hyp}) for $\ItlS\leq \ItlF$ with $\matV_0 \vec v^{\,(0)}=0$ satisfies as well the following equation for $0<\ItlS\leq \ItlF$
  \begin{equation}\label{eq:lm.gen.expr}
  \sum_{\ItlT=0}^{\ItlS-1}(-1)^{\ItlF-1-\ItlT} \binom{\ItlF}{\ItlT} \matD^{\ItlF-\ItlT} [\ItlF \matU-\matV_\ItlT] \vec v^{\,(\ItlT)}= (-1)^{\ItlF-\ItlS} \binom{\ItlF}{\ItlS} \matD^{\ItlF-\ItlS} \,  [\ItlS \matU-\matV_\ItlS] \vec v^{\,(\ItlS)}.
  \end{equation}

  We prove it by induction. Assuming $\ItlS=1$, we have that
  \begin{align*}
  (-1)^{\ItlF-1} \matD^{\ItlF} [\ItlF \matU-\matV_0] \vec v^{\,(0)}
  &=(-1)^{\ItlF-1} \matD^{\ItlF} [\ItlF \matU] \vec v^{\,(0)}\\
  &=(-1)^{\ItlF-1} \ItlF \matD^{\ItlF-1} \matD \matU \vec v^{\,(0)}\\
  &=(-1)^{\ItlF-1} \binom{\ItlF}{1} \matD^{\ItlF-1} (\matU-\matV_1) \vec v^{\,(1)},
  \end{align*}
  so that the base of the induction holds.
  Now assuming (\ref{eq:lm.gen.expr}) valid for $\ItlS<\ItlF$ we have for $\ItlS+1$ that
   \begin{align*}
  \sum_{\ItlT=0}^{\ItlS}(&-1)^{\ItlF-1-\ItlT} \binom{\ItlF}{\ItlT} \matD^{\ItlF-\ItlT} [\ItlF \matU-\matV_\ItlT] \vec v^{\,(\ItlT)}\\
  &=(-1)^{\ItlF-\ItlS} \binom{\ItlF}{\ItlS} \matD^{\ItlF-\ItlS} \,  [\ItlS \matU-\matV_\ItlS] \vec v^{\,(\ItlS)}+(-1)^{\ItlF-1-\ItlS} \binom{\ItlF}{\ItlS} \matD^{\ItlF-\ItlS} [\ItlF \matU-\matV_\ItlS] \vec v^{\,(\ItlS)}\\
  &=(-1)^{\ItlF-1-\ItlS} \binom{\ItlF}{\ItlS} \matD^{\ItlF-\ItlS} (\ItlF-\ItlS)  \matU \vec v^{\,(\ItlS)}\\
  &=(-1)^{\ItlF-1-\ItlS} \binom{\ItlF}{\ItlS} \matD^{\ItlF-\ItlS-1} (\ItlF-\ItlS) \matD \matU \vec v^{\,(\ItlS)}\\
  \intertext{and by equation (\ref{eq:lem.hyp}),}
  &=(-1)^{\ItlF-1-\ItlS} \binom{\ItlF}{\ItlS} \frac{\ItlF-\ItlS}{\ItlS+1}  \matD^{\ItlF-\ItlS-1}((\ItlS+1) \matU-\matV_{\ItlS+1}) \vec v^{\,(\ItlS+1)}\\
  &=(-1)^{\ItlF-(\ItlS+1)} \binom{\ItlF}{\ItlS+1}  \matD^{\ItlF-(\ItlS+1)}((\ItlS+1) \matU-\matV_{\ItlS+1}) \vec v^{\,(\ItlS+1)}
  \end{align*}
\end{proof}

The proof of next result follows closely the one of Lemma B.1 in \cite{seneta:1981}.

\begin{lem}\label{lm:invertibility}
The matrix $\matB_\ImF$, $\ImF>0$ in equation (\ref{eq:main.result.imp}) is invertible.
\end{lem}
\begin{proof}
By noticing that all the diagonal entries of the matrix $\diag\left(\tau_\IsF^{-1}(\ImF\mu_\IsF)\right)$ are strictly positive, in order to prove the invertibility of the matrix $\matB_\ImF$ we are left with proving the non-singularity of the matrix $(\matI-\matD_\ImF)$, with $$\matD_\ImF\defeq\diag\left(\tau_\IsF(\ImF\mu_\IsF)\right) \, \matQ.$$
We assume with no loss of generality that the states of the random environment are ordered according to the increasing values of $\beta_\IsF$, so that the states with null betas have the lowest indexes. Being $\beta>0$ we know that they are in number $K_0<K$.
If we compute the $\ImT$ power of the matrix $\matD_\ImF$ we get that in the limit it converges elementwise to
\begin{equation}\label{eq.diff.power.decomp}
\matD_\ImF^\ImT\to\left(
\begin{array}{cc}
  \vec1\,\transp{\vec\pi_0} & 0 \\
  0 & 0
\end{array}\right)
, \quad \textrm{as }\ImT\to\infty,
\end{equation}
where $\vec\pi_0$ is the vector containing only the first $K_0$ coordinates of the vector $\vec\pi$ and $\vec1$ is defined as a vector with all coordinates equal to 1 and whose dimension depends on the context.

The above result comes from the fact that $\matQ^\ImT\to\vec1\,\transp{\vec\pi}$ as $\ImT\to\infty$ and from knowing that $\tau_\IsF(\ImF\mu_\IsF)=1$ when  $\IsF\leq K_0$ and $\tau_\IsF^\ImT(\ImF\mu_\IsF)\to0$ when $K_0<\IsF\leq K$ as $\ImT\to\infty$.
For any $l>0$ the following relation is valid
$$(\matI-\matD_\ImF) (\matI+\matD_\ImF+\cdots+\matD_\ImF^{\ImT-1})=\matI-\matD_\ImF^\ImT,$$
and the matrix on the right side converges to
$$\matI-\matD_\ImF^\ImT\to\left(
\begin{array}{cc}
  \matI-\vec1\,\transp{\vec\pi_0} & 0 \\
  0 & \matI
\end{array}\right) \quad \textrm{as }\ImT\to\infty.$$
The determinant of the limit matrix in the right side of last relation is equal to the determinant of the matrix $(\matI-\vec1\,\transp{\vec\pi_0})$ that is positive by applying Lemma B.1 in \cite{seneta:1981} to the strictly substochastic matrix $\vec1\,\transp{\vec\pi_0}$.

Following the reasoning in Lemma B.1 in \cite{seneta:1981}, as the determinant is a continuous function in the space of matrices with elementwise convergence we get that there exists some $\ImT>0$ such that the determinant of the matrix $(\matI-\matD_\ImF^\ImT)$ is positive. Therefore it follows that the product of the determinants of the two matrix factors in the left side of equation (\ref{eq.diff.power.decomp}) has to be positive, that concludes the proof.
\end{proof}

\end{document}